\documentclass[11pt, leqno]{article}
\usepackage{amsmath, amsfonts}
\begin{document}

\author{Ajai Choudhry}
\title{A Note on the Quartic Diophantine Equation\\
$A^4+hB^4=C^4+hD^4$}
\date{}
\maketitle

\begin{abstract}
Integer solutions of the diophantine equation $A^4+hB^4=C^4+hD^4$ are known for all positive integer values of $h < 1000$. While a solution of the aforementioned diophantine equation for any arbitrary positive integer value of $h$ is not known, Gerardin and Piezas  found  solutions of this equation when $h$ is given by  polynomials of degrees 5 and 2 respectively. In this paper, we present several new solutions of this equation  when $h$ is given by   polynomials of degrees $2,\;3$  and 4.
\end{abstract}

This paper is concerned with the quartic diophantine equation,
\begin{equation}
A^4+hB^4=C^4+hD^4, \label{quarth}
\end{equation}
where $h$ is a given integer. Numerical solutions of Eq.~\eqref{quarth} for various values of $h$ were recorded by several authors (\cite[pp. 647-48]{Di}. A table of solutions for 75  integer values of $h \leq 101$ was compiled by Choudhry \cite{Ch}. This table was first  extended by Piezas \cite{Pi} to include all positive integer values of $h \leq 101$, and then by Tomita \cite{To} for all positive integer values of $h < 1000$ except $h=967$.  The missing solution for $h=967$ was supplied by Bremner (as mentioned by Tomita) and now solutions of  \eqref{quarth} are known for all positive integers $h < 1000.$ These numerical solutions confirm Choudhry's remark \cite{Ch} that Eq.~\eqref{quarth} seems to have a solution in integers for all positive integer values of $h$.  

There is no known method that would yield a solution in integers of \eqref{quarth} for any arbitrary value of $h$. When $h=2p^3(p^2-1)$, Gerardin had noted the solution $A=2p^2,\;B=p-1,\;C=2p,\;D=p+1$ (as quoted by Dickson \cite[p. 647]{Di}). Piezas has made a remark \cite{Pi2} from which it immediately follows that a solution of Eq.~\eqref{quarth} when $h=p^2-3$ is given by
\[ A=p^3+2p^2-3p-2,\;\;B=p^3-p- 2,\;\;C=p^3-2p^2-3p+2,\;\;D=p^3-p+2.\]
Further, Tomita \cite{To} has noted that when $h=p^4+q^4$, an obvious  solution of Eq.~\eqref{quarth} is given by $A=p^2,\,B=q,\,C=q^2,\,D=p$. 

In Table~\ref{tab1}  we present several new solutions of \eqref{quarth}  in which  $h$ is given by simple  polynomials of degrees $2,\;3$ and 4.  These solutions were obtained by experimenting with Eq.~\eqref{quarth} and may be readily verified by hand or by using any symbolic computation software such as MAPLE or Mathematica. We also note that when one nontrivial solution of \eqref{quarth} has been found, infinitely many integer solutions may be obtained following a procedure described in \cite{Ch}.

\begin{table}[tbh]
\caption{\bf Solutions of the equation $A^4+hB^4=C^4+hD^4$}
\label{tab1}
\begin{center}
\begin{tabular}{|c|c|c|c|c|} \hline
$h$&$A$&$B$&$C$&$D$\\
\hline
&&&&\\
$p^2+2  $ & $ p^3+2p+1  $ & $ p^2-p+1  $ & $ p^3+2p-1  $ & $ p^2+p+1   $ \\
$ p(p^2+4) $  &   $ p-2 $  &   $  2 $  &   $ p+2 $  &   $ 0 $\\$ 
 8p(p^2+1) $  &   $  p-1 $  &   $  1 $  &   $  p+1 $  &   $  0 $\\$
p^4-1 $  &   $ p $  &   $  0 $  &   $  1 $  &   $   1 $\\$ 
 2p^4-2 $  &   $  p^2+2p-1 $  &   $  p-1 $  &   $  p^2-2p-1 $  &   $ p+1 $\\$
 p^4+3p^2+1 $  &   $  p^2+p+1 $  &   $  p-1 $  &   $  p^2-p+1 $  &   $ p+1 $\\
\hline
\end{tabular}
\end{center}
\end{table}

We note that if a solution of Eq.~\eqref{quarth} when $h=\phi(p)$ is given by $A=A(p),\;B=B(p),\;C=C(p),\;D=D(p)$, then a rational solution of 
\eqref{quarth} when $h=q^4\phi(p/q)$ is given by $A=qA(p/q),\;B=B(p/q),\;C=qC(p/q),\;D=D(p/q)$, and integer solutions may be obtained by multiplying through by a constant. Thus, from Table~\ref{tab1}, we readily obtain solutions of \eqref{quarth} when the value of $h$ is given 
by any of the polynomials $(p^2+2q^2)q^2,\;pq(p^2+4q^2),\;8pq(p^2+q^2),\;p^4-q^4,\;2(p^4-q^4)$ and $p^4+3p^2q^2+q^4$.

\begin{center}
\Large
Acknowledgment
\end{center}
 
I wish to  thank the Harish-Chandra Research Institute, Allahabad for providing me with all necessary facilities that have helped me to pursue my research work in mathematics.

\bigskip

\noindent Postal Address: 13/4 A Clay Square, \\
\hspace{1.1in} Lucknow - 226001,\\
\hspace{1.1in}  INDIA \\
\medskip
\noindent e-mail address: ajaic203@yahoo.com

\end{document}